\documentclass[sts]{imsart}

\RequirePackage{amsthm,amsmath,amsfonts,amssymb}
\RequirePackage[numbers]{natbib}

\usepackage{verbatim}
\usepackage{comment}
\usepackage{graphicx}
\usepackage{enumerate}
\usepackage{color}
\usepackage{tikz}
\usepackage{hyperref}
\usepackage{mathrsfs}
\usepackage{euscript}
\usepackage{siunitx}

\usepackage{subcaption}
\usepackage{caption} 
\usetikzlibrary{shapes.geometric}
\usetikzlibrary{backgrounds}

\startlocaldefs

\DeclareGraphicsExtensions{.pdf,.jpg,.png}

\theoremstyle{plain}
\newtheorem{theorem}{Theorem}

\newtheorem{conjecture}[theorem]{Conjecture}

\theoremstyle{definition}

\newcommand{\G}{\mathcal G}
\newcommand{\const}{\mathscr C}
\newcommand{\dd}{\,\mathrm{d}}

\newcommand{\T}{\mathsf{T}}
\newcommand{\E}{\mathcal E}
\newcommand{\V}{\mathcal V}
\newcommand{\CfourA}{\begin{tikzpicture}[thick, scale=0.8, every node/.style={circle, minimum size = 3pt, scale=0.8}]
		\node[shape=circle,draw=black] (1) at (0,2) {$SL$};
		\node[shape=circle,draw=black] (2) at (2,2) {$SW$};
		\node[shape=circle,draw=black] (3) at (0,0) {$PL$};
		\node[shape=circle,draw=black] (4) at (2,0) {$PW$};
		
		\draw (2) -- (4);
		\draw (3) -- (4);
		\draw (1) -- (3);
		\draw (2) -- (1);
		
\end{tikzpicture}}
\newcommand{\CfourB}{\begin{tikzpicture}[thick, scale=0.8, every node/.style={circle, minimum size = 3pt, scale=0.8}]
		\node[shape=circle,draw=black] (1) at (0,2) {$SL$};
		\node[shape=circle,draw=black] (2) at (2,2) {$SW$};
		\node[shape=circle,draw=black] (3) at (0,0) {$PL$};
		\node[shape=circle,draw=black] (4) at (2,0) {$PW$};
		
		\draw (1) -- (4);
		\draw (1) -- (3);
		\draw (2) -- (3);
		\draw (2) -- (4);
\end{tikzpicture}}
\newcommand{\CfourC}{\begin{tikzpicture}[thick, scale=0.8, every node/.style={circle, minimum size = 3pt, scale=0.8}]
		\node[shape=circle,draw=black] (1) at (0,2) {$SL$};
		\node[shape=circle,draw=black] (2) at (2,2) {$SW$};
		\node[shape=circle,draw=black] (3) at (0,0) {$PL$};
		\node[shape=circle,draw=black] (4) at (2,0) {$PW$};
		
		\draw (1) -- (4);
		\draw (3) -- (4);
		\draw (2) -- (3);
		\draw (2) -- (1);
\end{tikzpicture}}

\DeclareMathOperator{\Iss}{Iss}
\DeclareMathOperator{\tr}{tr}
\DeclareMathOperator{\diag}{diag}

\newcommand{\rev}[1]{#1}

\endlocaldefs

\makeatletter
\def\journal@name{}
\makeatother

\begin{document}

\begin{frontmatter}
\title{On a conjecture of Roverato regarding\\ $G$-Wishart normalising constants}
\runtitle{On Roverato's $G$-Wishart conjecture}

\begin{aug}
\author[A]{\fnms{Ching}~\snm{Wong}\ead[label=e1]{ching.wong@unibas.ch}},
\author[B]{\fnms{Giusi}~\snm{Moffa}\ead[label=e2]{giusi.moffa@unibas.ch}}
\and
\author[C]{\fnms{Jack}~\snm{Kuipers}\ead[label=e3]{jack.kuipers@bsse.ethz.ch}}
\address[A]{Department of Mathematics and Computer Science, University of Basel, Switzerland\printead[presep={\ }]{e1}.}
\address[B]{Department of Mathematics and Computer Science, University of Basel, Switzerland\printead[presep={\ }]{e2}.}
\address[C]{Department of Biosystems Science and Engineering, ETH Zurich, Switzerland\printead[presep={\ }]{e3}.}
\end{aug}

\begin{abstract}
The evaluation of $\G$-Wishart normalising constants is a core component for Bayesian analyses for Gaussian graphical models, but remains a computationally intensive task in general. Based on empirical evidence, Roverato [Scandinavian Journal of Statistics, 29:391--411 (2002)] observed and conjectured that such constants can be simplified and rewritten in terms of constants with an identity scale matrix. In this note, we disprove this conjecture for general graphs by showing that the conjecture instead implies an independently-derived approximation for certain ratios of normalising constants. \rev{We further show that the conjecture is actually a saddle-point approximation, allowing us to derive the next-order correction which enables more accurate evaluation.}
\end{abstract}


\end{frontmatter}
\section{Introduction}
	
Since the $\G$-Wishart distribution acts as the conjugate prior for Gaussian graphical models, computing their posterior for Bayesian analyses involves evaluating the marginal likelihood related to the \emph{$\G$-Wishart normalising constant} \cite{Roverato}:
\begin{equation} \label{eq:C_G}
	\const_\G(\delta, D) := \int_{\mathbb S^n_{++}(\G)} \det(K)^{\frac{ \delta-2}{2}} e^{-\frac{\tr(KD)}2} \dd K,
\end{equation}
for $\delta > 0$ and $D \in \mathbb S^n_{++}$ where the set $\mathbb S^n_{++}$ denotes the cone of all $n$ by $n$ real symmetric positive definite matrices. For a graph $\G$ with vertices $\V(\G) = \{v_1, \ldots, v_n\}$, the set $\mathbb S^n_{++}(\G)$ contains all the matrices in $\mathbb S^n_{++}$ whose $(\mu,\nu)$-entry is zero whenever there is no edge connecting $v_\mu$ and $v_\nu$ in $\G$.
	
For general $D$, the constant $\const_\G(\delta, D)$ is only known for \emph{chordal graphs}, i.e., graphs whose induced cycles have length 3 \cite{Dawid}, and graphs that can be made chordal with the addition of a single edge \cite{GWishart}. Otherwise, directly estimating the constant using existing Monte Carlo methods \cite{Roverato, Atay-Kayis, mw19} requires sampling from the $\G$-Wishart distribution and can be computationally intensive for larger networks. \rev{More recently, a Monte Carlo approach was developed using Fourier analysis, which, while still intensive, can be orders of magnitude more efficient \cite{GWishart}.}

For Bayesian sampling of Gaussian graphical models, one popular approach is to employ Markov chain Monte Carlo (MCMC) to target the posterior distribution of graphs, with respect to a given dataset. To avoid the expensive computation of the marginal likelihoods, a widely used method is to instead sample in the joint distribution of the graph $\G$ and precision matrix $K$ and update the precision matrix in each MCMC iteration \cite{mml23}. The former task still involves the computation of the normalising constants when $D$ is the identity matrix. This however is generally more efficient and exact formulae are known for some graphs \cite{Uhler, GWishart}. To avoid directly computing the normalising constant even for the identity matrix, MCMC schemes then use exchange or auxiliary variable approaches to estimate and remove all normalising constants from the acceptance ratios \cite{mw15, Willem}. This still requires sampling from the $\G$-Wishart distribution \cite{lenkoski2013direct}.

Rather than targeting the full posterior, approximations have been developed for the ratio of normalising constants \cite{mml23} based on the Cholesky decomposition of the $\G$-Wishart distribution. In particular, let $\G^*$ be a graph on $n$ vertices and let $\G$ be the graph obtained from $\G^*$ by removing one of its edges $e$, then 
\begin{equation} 
	\const_{\G} (\delta, I_n) \approx \dfrac{\const_{\G^*}(\delta, I_n)}{2\pi^{\frac12}} \dfrac{ \Gamma\left(\frac{\delta+s}{2}\right)}{\Gamma\left(\frac{\delta+s+1}{2}\right)}, \label{eq:ratio}
\end{equation}
where $s$ is the number of common neighbours of the end vertices of $e$ in $\G^*$. While exact for some graphs, the approximation does not hold in general, although it can be quite good \cite{mml23}.

Even with approximations or exchange algorithms, sampling in the joint space of $(\G ,K)$ and updating the precision matrix adds considerable computational overhead. It would therefore be desirable to derive a formula for $\const_\G(\delta, D)$ based on $\const_\G(\delta, I_n)$, and in 2002, Roverato \cite{Roverato} studied the function
\begin{equation}
    h_\G(\delta, D) := 2^{\frac n 2} \det(\Iss_\G(D^\G))^{-\frac12} \det(D^\G)^{-\frac{\delta-2}{2}} \const_\G(\delta, D)^{-1}, \nonumber
\end{equation}
where $D^\G$ is the \emph{PD-completion of $D$ with respect to $\G$}, and $\Iss_\G(D^\G)$ is the \emph{Isserlis matrix} of $D^\G$ with respect to $\G$ (definitions in Section~\ref{sec:conj}). PD-completion can be performed by iteratively updating $D$ \cite{speed1986gaussian}, so that both determinants can be computed relatively easily. Based on some \rev{simulations}, Roverato observed \rev{\citep[p.~403]{Roverato} ``these tables provide an empirical evidence that, as well as for
the decomposable case [chordal graphs], also when G is non-decomposable the normalizing constant [$h_\G(\delta, D)$] does not
depend on $D$''.}  Together with the fact that $\det(\Iss_\G(I_n)) = 2^n$, this leads to a very desirable conjecture \rev{\citep[see also][p.~400 and p.~408]{Roverato}}:
\begin{conjecture}[Roverato] \label{conj}
	For graph $\G$ with $n$ vertices, $D \in \mathbb S^n_{++}$ and real number $\delta > 0$,
	\begin{equation}
		\const_\G(\delta, D) = 2^{\frac{n}{2}} \det(\Iss_\G(D^\G))^{-\frac12} \det(D^\G)^{-\frac{\delta-2}2} \const_\G(\delta, I_n). \nonumber
	\end{equation}
\end{conjecture}

Recently, we devised a new way of evaluating $\const_\G(\delta, D)$ via Fourier analysis \cite{GWishart}. \rev{The Fourier analysis turns the integral of Equation~(\ref{eq:C_G}) over the edges \emph{present} in $\G$ into the dual integral over the \emph{absent} edges:
\begin{equation} \label{eq:I_G}
\const_\G(\delta, D) = C^{\gamma}_n(\G) \int_{\mathbb S^{n}(\G^c)} \det(D+iT)^{-\gamma} \dd T,
\end{equation}
for $\gamma = \frac{n+\delta-1}{2}$, where $C^{\gamma}_n(\G)$ is a constant depending on $n$, $\gamma$ and the number of edges in $\G$, and where $\mathbb S^{n}(\G^c)$ is the set of $n$ by $n$ real symmetric matrices $T$ whose entries are all zero except where there are edges in the complement $\G^c$ (so where edges are absent in $\G$). For completeness, we sketch the derivation in Appendix A.}

Using this approach, we \rev{first} show in this note that Roverato's conjecture implies the approximation in Equation~(\ref{eq:ratio}) for all non-empty graphs $\G^*$, and the conjecture is therefore not true in general. \rev{We then show that the conjecture corresponds to an asymptotic approximation and derive the next-order correction.}
	
	\section{An equivalent form of the conjecture} \label{sec:conj}
	
	We first recall some definitions from \cite{Roverato}. Let $D \in \mathbb S^n_{++}$ and let $\G$ be a graph with vertex set $\V(\G) = \{v_1, \ldots, v_n\}$ and edge set $\E(\G) \subseteq \{\{v_\mu, v_\nu\}: \mu \neq \nu\}$. Let $m = |\E(\G)|$.
	
	The \emph{PD-completion of $D$ with respect to $\G$}, denoted by $D^\G$, is the unique matrix in $\mathbb S^n_{++}$ such that $(D^\G)^{-1} \in \mathbb S^n_{++}(\G)$ and the matrices $D$ and $D^\G$ agree on the diagonal entries	as well as the $(\mu, \nu)$-entries whenever $\{v_\mu, v_\nu\}$ is an edge in $\G$. The existence and uniqueness of $D^\G$ is established in \cite{Grone} while iterative schemes to compute it are detailed in \cite{speed1986gaussian}.
	
	The \emph{Isserlis matrix of $D$ with respect to $\G$}, denoted by $\Iss_\G(D)$, is the symmetric matrix in $\mathbb R^{(n+m) \times (n+m)}$ indexed by $W \times W$, where
	\begin{align}
	W = & \{(\mu, \mu) : 1 \leq \mu \leq n\} \cup \nonumber \\ & 
	\{(\mu, \nu) : \mu < \nu \text{ and } \{v_\mu, v_\nu\} \in \E(\G)\}, \nonumber
	\end{align}
	whose entry associated with the index $((\mu, \nu), (\mu', \nu')) \in W \times W$ is $d_{\mu, \mu'} d_{\nu, \nu'} + d_{\mu\nu'} d_{\mu'\nu}$. If $\G$ is the \emph{complete graph} (i.e., $\G$ has $\binom n 2$ edges), then we write $\Iss(D)$ in place of $\Iss_\G(D)$.
	
	The following identities regarding the Isserlis matrices are well known:
	\begin{itemize}
		\item $\Iss(D)$ is invertible and \\
		 $(\Iss(D))^{-1} = \Iss(I_n)^{-1} \Iss(D^{-1}) \Iss(I_n)^{-1}$ \cite{RoveratoWhittaker},
		\item $\Iss(I_n) = \diag(\underbrace{2, \ldots, 2}_{n}, \underbrace{1, \ldots, 1}_{\binom n 2})$,
		\item $\det(\Iss(D)) = 2^n \det(D)^{n+1}$.
	\end{itemize}
	
	Since we are only interested in the determinants of the Isserlis matrices in this note, the ordering of the rows and columns does not matter. Thus, we may assume that the top left $n$ by $n$ block of $\Iss_\G(D)$ corresponds to the first part of $W$, i.e., $\{(\mu, \mu) : 1 \leq \mu \leq n\}$. Given a graph $\G$ on $n$ vertices and $m$ edges, we may assume that the top left $n+m$ by $n+m$ block of $\Iss(D)$ coincides with $\Iss_\G(D)$.
	
	For $1\leq \mu \leq \nu \leq n$, we use $D[\mu, \nu]$ to denote the $\nu-\mu+1$ by $\nu-\mu+1$ submatrix of $D$ formed by selecting the rows and columns ranging from the $\mu$-th to the $\nu$-th indices. 
	
	We consider the matrix $\Iss(D^\G)$ as a 2 by 2 block matrix, whose top left block is $\Iss_\G(D^\G)$. Then, the determinant of $\Iss(D^\G)$ is
	\begin{align}
		& 2^n \det(D^\G)^{n+1} = \det(\Iss(D^\G)) \nonumber \\
		& = \det(\Iss_\G(D^\G)) \nonumber \\
		& \qquad \times \det\left((\Iss(D^\G))^{-1}\left[n+m+1,n+\binom n 2\right]\right)^{-1} \nonumber \\
		& = \det(\Iss_\G(D^\G)) \nonumber \\
		& \qquad \times \det\left(\Iss((D^\G)^{-1})\left[n+m+1,n+\binom n 2\right]\right)^{-1}. \nonumber
	\end{align}
    where the first step follows from the block determinant involving the Schur complement, which is expressed in terms of a block of the inverse, while the second step follows from the identity above concerning the inverse of the Isserlis matrix. 
	
	Therefore, Roverato's conjecture is equivalent to
	\begin{align}
		\const_\G(\delta, D) = & \det\left(\Iss((D^\G)^{-1})\left[n+m+1,n+\binom n 2\right]\right)^{-\frac12} \nonumber \\
		& \times \det(D^\G)^{-\frac{\delta-2}2-\frac{n+1}{2}} \const_\G(\delta, I_n), \label{eq:conj}
	\end{align}
	where $\G$ is a graph with $n$ vertices and $m$ edges, $D \in \mathbb S^n_{++}$, $\delta > 0$ is a real number. \rev{We note that this formulation involves the Isserlis matrix of the edges absent in the graph $\G$ so that it aligns with the Fourier transformed integral in Equation~(\ref{eq:I_G}).}
	
	\section{The conjecture implies the ratio estimate}
	
	In this section, we show that Roverato's conjecture from Equation~(\ref{eq:conj}) implies the ratio approximation in Equation~(\ref{eq:ratio}).
	
	Let $\G^*$ be a graph with vertex set $\V(\G^*) = \{v_1, \ldots, v_n\}$ such that $\{v_1, v_2\} \in \E(\G^*)$. Let $\G$ be the graph resulting from the removal of edge $\{v_1, v_2\}$ from $\G^*$. Let $\delta > 0$ be a real number and $D \in \mathbb S^n_{++}$. Let $m = |\E(\G)|$.
	
	The main result in \cite{GWishart} implies that
	\begin{eqnarray} \label{eq:missing_edge}
		\const_{\G} (\delta, I_n) = \dfrac{1}{2\pi} \int_{\mathbb S} \const_{\G^*}(\delta, I_n+iT) \dd T, 
	\end{eqnarray}
	where $\mathbb S$ is the set of $n$ by $n$ matrices $T$ whose entries are all zero except $t_{1,2} = t_{2,1} \in \mathbb R$. The normalising constant $\const_{\G^*}(\delta, I_n+iT)$ is defined using analytic continuation. \rev{Compared to Equation~(\ref{eq:I_G}) where we integrate over all the absent edges in $\G$ (relative to a complete graph), here we just integrate over the one edge absent from $\G^*$.}

\rev{We assume an analytic continuation of Roverato's conjecture, consistent with the analytic continuation used to define $\const_{\G^*}(\delta, I_n+iT)$ in the integrand of Equation~(\ref{eq:missing_edge}), for which we use Equation~(\ref{eq:conj}) where the matrix argument is the non-identity $I_n+iT$}:
	\begin{align} \label{eq:conj2}
		& \const_{\G} (\delta, I_n) = \dfrac{\const_{\G^*}(\delta, I_n)}{2\pi} \int_{\mathbb S} \\
		& \det\left(\Iss(((I_n+iT)^{\G^*})^{-1})\left[n+m+2,n+\binom n 2\right]\right)^{-\frac12}  \nonumber \\
		& \qquad \times \det((I_n+iT)^{\G^*})^{-\frac{\delta-2}2-\frac{n+1}{2}}\dd T. \nonumber 
	\end{align}
	Since
	\begin{equation}
		(I_n + iT)^{\G^*} = I_n + iT \nonumber
	\end{equation}
	and
	\begin{align}
		& (I_n + iT)^{-1} = \diag\left(\frac{1}{1+t^2}\begin{pmatrix}
			1 & -it \\
			-it & 1
		\end{pmatrix}, I_{n-2}\right), \nonumber
	\end{align}
    the bottom right block of the matrix
	\begin{align}
		& \Iss(((I_n+iT)^{\G^*})^{-1}) = \nonumber \\
		 & \Iss\left(\diag\left(\frac{1}{1+t^2}\begin{pmatrix}
			1 & -it \\
			-it & 1
		\end{pmatrix}, I_{n-2}\right)\right) , \nonumber
	\end{align}
	contains those entries associated with the non-adjacent pairs of vertices in $\G^*$. As the ordering of the pairs does not matter, we may assume that the first $2w$ pairs to be
	\begin{equation}
	(v_1, v_{a_1}), (v_2, v_{a_1}), 
	\ldots, (v_1, v_{a_w}), (v_2, v_{a_w}), \nonumber
	\end{equation}
	where $w$ is the number of vertices (excluding $v_1, v_2$) that are adjacent to neither $v_1$ nor $v_2$ in $\G$. Next, we have $(v_1, v_{b_1}), \ldots, (v_1, v_{b_x})$, where $x$ is the number of vertices (excluding $v_2$) that are not adjacent to $v_1$, but adjacent to $v_2$ in $\G$. It is then followed by $(v_2, v_{c_1}), \ldots, (v_2, v_{c_y})$ similarly. We remark that $s := n-2-(w+x+y)$ is the number of common neighbours of $v_1, v_2$ in $\G$. The final pairs are those not involving the vertices $v_1$ or $v_2$.
	
	The matrix block of interest is
	\begin{align}
		\diag\Bigg(\underbrace{\frac{1}{1+t^2}\begin{pmatrix}
				1 & -it \\
				-it & 1
			\end{pmatrix}, \ldots, \frac{1}{1+t^2}\begin{pmatrix}
				1 & -it \\
				-it & 1
		\end{pmatrix}}_w,  \nonumber \\
		\frac{1}{1+t^2}I_{x+y}, I_{\binom n 2 - m - 1 - 2w - x - y} \Bigg), \nonumber
	\end{align}
	whose determinant is
	\begin{align}
		&\det\left(\Iss(((I_n+iT)^{\G^*})^{-1})\left[n+m+2,n+\binom n 2\right]\right) \nonumber \\
		& = (1+t^2)^{-w-x-y}. \nonumber
	\end{align}
	Therefore, it follows from Equation~(\ref{eq:conj2}) that
	\begin{eqnarray}
		\const_{\G} (\delta, I_n) & = & \dfrac{\const_{\G^*}(\delta, I_n)}{2\pi} \int_{\mathbb R} (1+t^2)^{\frac{w+x+y}{2}-\frac{\delta-2}2-\frac{n+1}{2}} \dd t \nonumber \\
		& = & \dfrac{\const_{\G^*}(\delta, I_n)}{2\pi} \int_{\mathbb R} (1+t^2)^{-\frac{\delta+s+1}{2}}  \dd t \nonumber \\
		& = & \dfrac{\const_{\G^*}(\delta, I_n)}{2\pi^{\frac12}} \dfrac{\Gamma\left(\frac{\delta+s}{2}\right)}{\Gamma\left(\frac{\delta+s+1}{2}\right)}, \nonumber
	\end{eqnarray}
	which is the same as Equation~(\ref{eq:ratio}).
	
	In summary, we prove in this section that Conjecture~\ref{conj} for a non-empty graph $\G^*$ implies the ratio in Equation~(\ref{eq:ratio}) for all graphs $\G$ which have one fewer edge than $\G^*$. \rev{For the proof, we relied on the conjecture from Equation~(\ref{eq:conj}) in terms of the absent edges to match the integral in Equation~(\ref{eq:missing_edge}) also over the absent edge.}

\begin{figure}
	\includegraphics[width=0.48\textwidth]{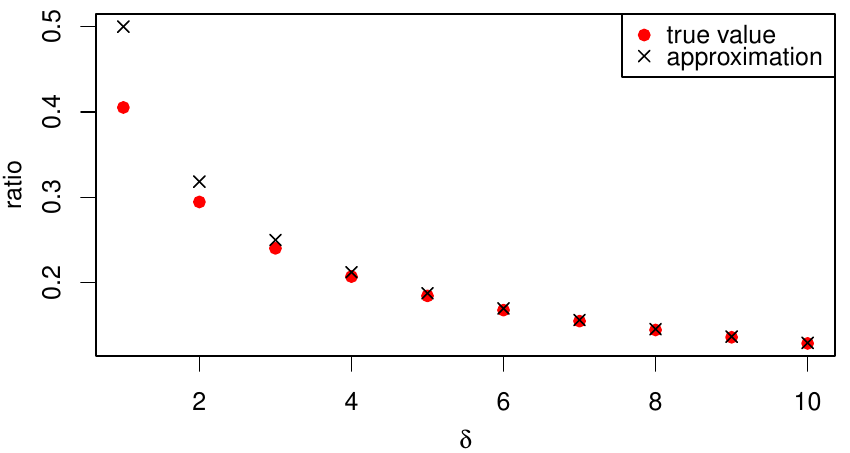}
		\caption{For the first ten positive integers $\delta$, we compare the correct ratio from Equation~(\ref{eq:true}) and the approximation of Equation~(\ref{eq:approx}).} \label{fig:ratio} 
\end{figure}

\section{A counter example}

    \begin{figure*}[t]
		\begin{subfigure}{1\textwidth}
			\centering
			\hspace{0.1\textwidth} \CfourA \hfill \CfourB \hfill \CfourC \hspace{0.1\textwidth}
			\caption{}
			\label{fig:C4a}
		\end{subfigure}
		\\
		\begin{subfigure}{1\textwidth}
			\centering
			\includegraphics[width=0.32\textwidth]{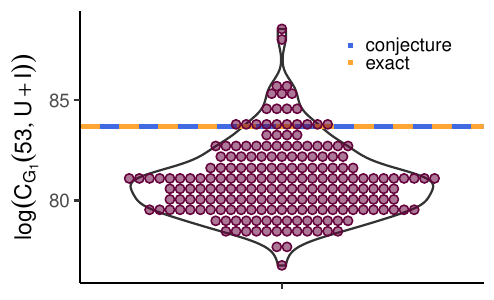} \includegraphics[width=0.32\textwidth]{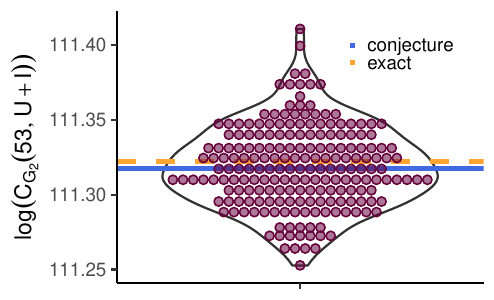} \includegraphics[width=0.32\textwidth]{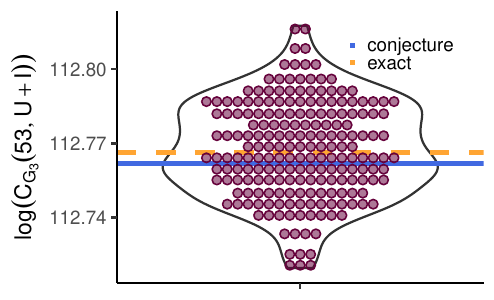}
			\caption{}
			\label{fig:C4b}
		\end{subfigure}

 		\caption{(a) The graphs $\G_1$, $\G_2$, $\G_3$ (from left to right), are the non-chordal graphs for Fisher's Iris Virginica dataset. (b) Violin plots of the estimates of the values of $\log(\const_{\G_j}(53, U+I_4))$, where $j = 1,2,3$, using Monte Carlo integration \cite{Atay-Kayis,mw19} with $1000$ samples and for 200 different seeds. The horizontal lines represent the values obtained using Roverato's conjecture and the exact one-dimensional integral in \cite{GWishart}.} \label{fig:C4}
	\end{figure*}

Take $\G^*$ be the cycle of length 4 and $\G$ be the path of length 3. For both graphs, exact formulae for the normalising constants (when $D$ is the identity matrix) are known \cite{Uhler}:
\begin{eqnarray}
	\const_{\G}(\delta, I_4) & = & 2^{2\delta+3} \pi^{\frac32} \Gamma\left(\frac{\delta+1}{2}\right)^3 \Gamma\left(\frac{\delta}{2}\right), \nonumber \\
	\const_{\G^*}(\delta, I_4) & = & 2^{2\delta+4} \dfrac{\pi^2 \Gamma\left(\frac{\delta}{2}\right) \Gamma\left(\frac{\delta+1}{2}\right) \Gamma\left(\frac{\delta+2}{2}\right)^3}{\Gamma\left(\frac{\delta+3}{2}\right)}. \nonumber
\end{eqnarray}
Therefore, the correct ratio of the normalising constants is
\begin{equation} 
	\dfrac{\const_{\G} (\delta, I_4)}{\const_{\G^*}(\delta, I_4)} = \dfrac{\Gamma\left(\frac{\delta+1}{2}\right)^2 \Gamma\left(\frac{\delta+3}{2}\right)}{2\pi^{\frac12}\Gamma\left(\frac{\delta+2}{2}\right)^3}, \label{eq:true}
\end{equation}
while the approximation from (\ref{eq:ratio}) is
\begin{equation} 
\dfrac{\const_{\G} (\delta, I_4)}{\const_{\G^*}(\delta, I_4)} \approx \dfrac{ \Gamma\left(\frac{\delta}{2}\right)}{2\pi^{\frac12} \Gamma\left(\frac{\delta+1}{2}\right)}. \label{eq:approx}
\end{equation}

    By Stirling's approximation, the relative error of the above approximation is
    \begin{align}
        & \dfrac{\Gamma\left(\frac{\delta+2}{2}\right)^3 \Gamma\left(\frac{\delta}{2}\right)}{\Gamma\left(\frac{\delta+1}{2}\right)^3\Gamma\left(\frac{\delta+3}{2}\right)} - 1 \nonumber \\
        & = \left(1+\dfrac{1}{4\delta}+\dfrac{1}{32\delta^2}\right)^3\left(1-\dfrac{3}{4\delta}+\dfrac{25}{32\delta^2}\right)\nonumber \\
        & \qquad \times (1+O(\delta^{-3}))^4 - 1 \nonumber \\
        & \sim \dfrac{1}{2\delta^2}, \quad \text{ as } \delta \to \infty. \nonumber
    \end{align}
    Therefore the error tends to zero as $\delta$ grows. Moreover, the first-order terms cancel so that the error decreases rapidly for increasing $\delta$. Figure~\ref{fig:ratio} compares the ratios for small values of $\delta$.

    \begin{table*}[t!]
		\caption{The estimated values of $\log(\const_{\G_j}(53, U+I_4))$, for $j = 1, 2, 3$. The graphs are displayed in Figure~\ref{fig:C4a}.}
		\begin{tabular}{l|S[table-format=3.4]|S[table-format=3.4]|S[table-format=3.4]}
			& {$\log(\const_{\G_1}(53, U+I_4))$} & {$\log(\const_{\G_2}(53, U+I_4))$} & {$\log(\const_{\G_3}(53, U+I_4))$}\\
			\hline
			{exact value \cite{GWishart}} & 83.6851 & 111.3223 & 112.7664 \\
			\hline
			{conjectured value} & 83.6836 & 111.3175 & 112.7618
		\end{tabular} \label{tab:logC}
	\end{table*}

    In the context of Bayesian inference, we need to compute the normalising constant $\const_\G(\delta + N, U+D)$, where $N$ is the number of data points, $U$ is the scatter matrix of the dataset, $\delta$ and $D$ are parameters of the prior distribution, often taken to be 3 and the identity matrix respectively. Since the parameter of the normalising constant in the posterior is $(\delta + N)$, with reasonable amounts of data this parameter will accordingly be quite large, and hence we may be in a regime where the approximation implied by Conjecture~\ref{conj} is quite accurate.

    As illustration, we use Fisher's Iris Virginica dataset, which contains 4 variables (SL, SW, PL, PW) and 50 data points. There are three non-chordal graphs, as shown in Figure~\ref{fig:C4a}. In \cite{GWishart}, we derived a new way to write the normalising constant $\const_{\G}(53, U+I)$ of these graphs as a one-dimensional integral, which can be computed numerically. For each of these three graphs, we compare this value and the conjectured value in Table~\ref{tab:logC}. Further, we also compare with values obtained by Monte Carlo integration \cite{Atay-Kayis} in Figure~\ref{fig:C4b}. It can be seen that Roverato's conjecture provides a very good estimate of the normalising constant, especially compared to the stochastic noise of the Monte Carlo integration.


\section{The conjecture as a saddle-point approximation}

\rev{To explore the approximation induced by Conjecture 1 in more detail, we return to the integral in Equation~(\ref{eq:I_G}) applied to the PD-completion of the matrix $D$
\begin{align} \label{eq:int_IDG}
&\int_{\mathbb S^{n}(\G^c)} \det(D^{\G}+iT)^{-\gamma} \dd T \\
& = \det(D^\G)^{-\gamma} \int_{\mathbb S^{n}(\G^c)} \det(I_n + i(D^{\G})^{-1}T)^{-\gamma} \dd T .\nonumber
\end{align}
For the integrand we reformulate 
\begin{align}
\det(I_n+i(D^{\G})^{-1}T)^{-\gamma} & = e^{-\gamma\log\det(I_n+i(D^{\G})^{-1}T)} \nonumber \\
= e^{-\gamma\tr\log(I_n+i(D^{\G})^{-1}T)} & \approx e^{-\frac{\gamma}{2}\tr((D^{\G})^{-1}T)^{2}} , \nonumber
 \end{align}
Since we used the PD-completion, the matrix $(D^{\G})^{-1}$ is zero exactly where $T$ is not (by construction) so that the linear terms vanish. This expansion is formal and corresponds to a Laplace approximation around $T=0$. We consider higher-order terms later and focus here on performing the saddle-point approximation by only keeping up to quadratic terms in the exponent.

Define $\boldsymbol{t}$ to be a vector corresponding to the non-zero element pairs of $T$, our quadratic form is
\begin{equation} \nonumber
\tr((D^{\G})^{-1}T)^{2} = 2\boldsymbol{t}^\T\Sigma^{-1}\boldsymbol{t},
\end{equation}
with
\begin{align} \label{eq:sigma_def}
\Sigma^{-1}_{ij} = & (D^{\G})^{-1}_{\E_{i1}, \E_{j1}}(D^{\G})^{-1}_{\E_{i2}, \E_{j2}}  \\ 
& + (D^{\G})^{-1}_{\E_{i1}, \E_{j2}}(D^{\G})^{-1}_{\E_{i2}, \E_{j1}} , \nonumber 
\end{align}
where $\E_{i1}$ and $\E_{i2}$ are the row and column position in $T$ corresponding to the $i$-th entry in $\boldsymbol{t}$. The precision matrix $\Sigma^{-1}$ is then exactly the Isserlis matrix of $(D^{\G})^{-1}$ for the absent edges represented in $T$, giving
\begin{align} \label{eq:iss_approx}
&\det\left(\Sigma^{-1}\right) \\
& = \det\left(\Iss((D^\G)^{-1})\left[n+m+1,n+\binom n 2\right]\right) . \nonumber
\end{align}
Combining these components
\begin{align}
&\int_{\mathbb S^{n}(\G^c)} \det(D^{\G}+iT)^{-\gamma} \dd T \nonumber\\
& \approx \det(D^\G)^{-\gamma} \int  e^{-\gamma \boldsymbol{t}^\T\Sigma^{-1}\boldsymbol{t}} \dd \boldsymbol{t} \nonumber \\
& = \frac{\det(D^\G)^{-\gamma}}{\det(\Sigma)^{-\frac12}} \int e^{-\gamma \boldsymbol{t}^\T \boldsymbol{t}} \dd \boldsymbol{t}, \nonumber
\end{align}
through a change of variables. Rather than actually perform the Gaussian integral, we run through the same approximations with the identity matrix
\begin{align}
\int_{\mathbb S^{n}(\G^c)} \det(I_n+iT)^{-\gamma} \dd T \approx \int e^{-\gamma \boldsymbol{t}^\T \boldsymbol{t}} \dd \boldsymbol{t}, \nonumber
\end{align}
to arrive at
\begin{align}
&\int_{\mathbb S^{n}(\G^c)} \det(D^{\G}+iT)^{-\gamma} \dd T \nonumber\\
& \approx \frac{\det(D^\G)^{-\gamma}}{\det(\Sigma)^{-\frac12}} \int_{\mathbb S^{n}(\G^c)} \det(I_n+iT)^{-\gamma} \dd T . \nonumber
\end{align}
Substituting for Equation~(\ref{eq:iss_approx}) and using Equation~(\ref{eq:I_G}) we recover the conjecture in the form of Equation~(\ref{eq:conj}).

\begin{figure*}
\begin{tabular}{cc}
\includegraphics[width=0.48\textwidth]{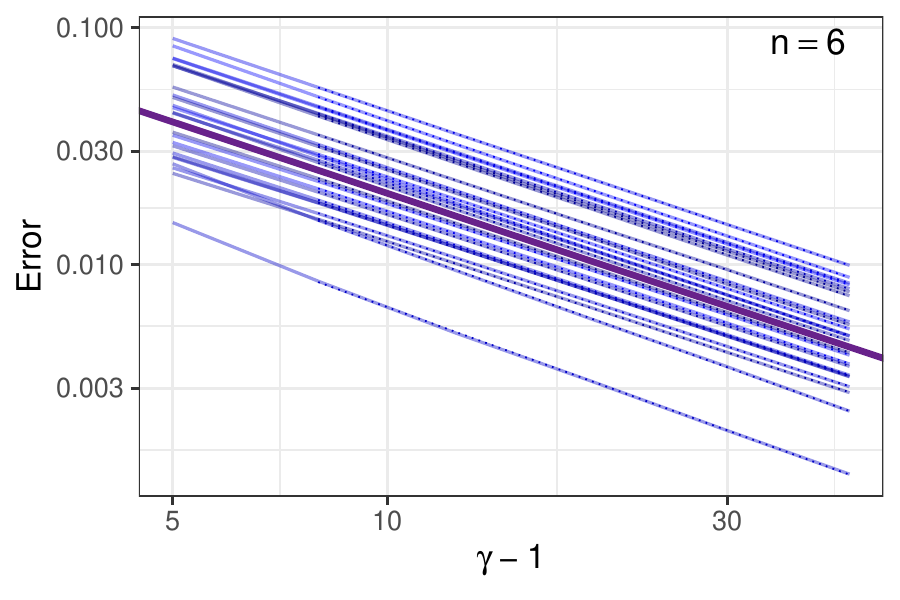} & \includegraphics[width=0.48\textwidth]{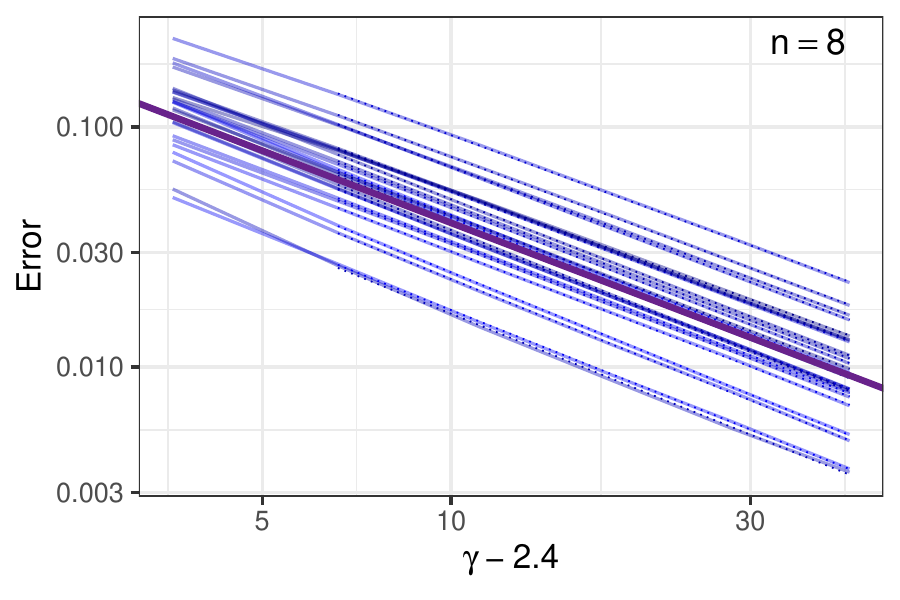}\\
\includegraphics[width=0.48\textwidth]{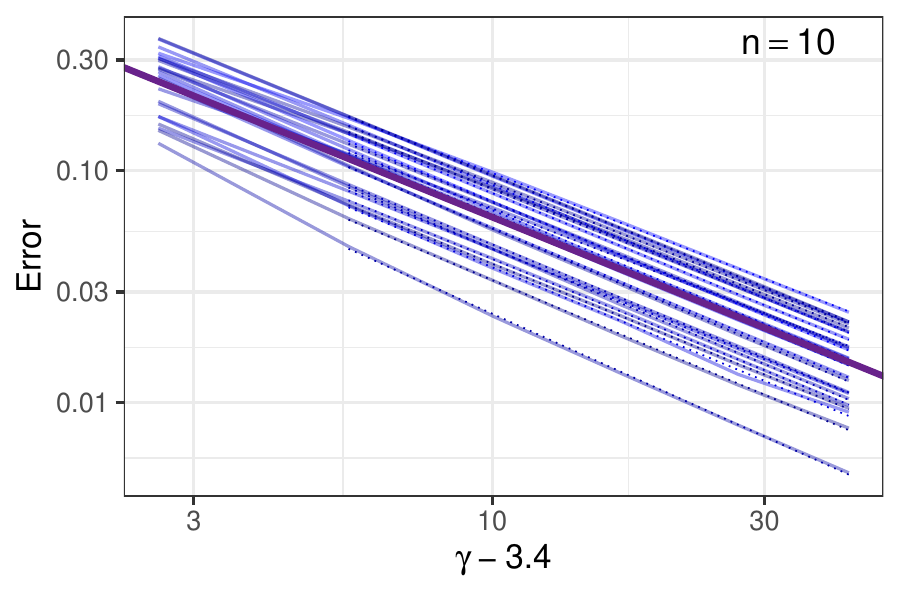} & \includegraphics[width=0.48\textwidth]{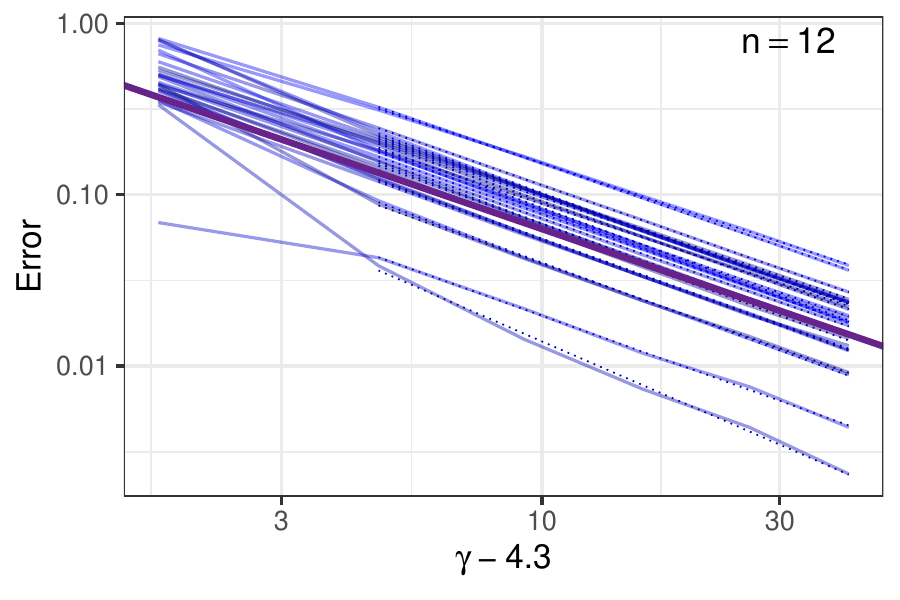}
\end{tabular}
\caption{\rev{For 25 different random prime graphs of size 6, 8, 10, and 12, with randomly generated $D$, we show the relative error from the approximation for increasing $\gamma$ in a log-log plot. Each line is for different graph, with the dotted lines derived from a linear fit. The thick purple line is an indicative line with a slope of $-1$. To improve the linear fits by mitigating higher-order corrections, we shifted $\gamma$ by $-1$, $-2.4$, $-3.4$, and $-4.3$, for the different $n$ respectively.}} \label{fig:convergeapprox} 
\end{figure*}

The conjecture therefore corresponds to a change of variables so that $\const_\G(\delta, D)$ matches $\const_\G(\delta, I_n)$ perfectly up to the second order terms in a saddle-point approximation. Errors then accumulate from the higher order terms where they no longer match. With the anti-symmetry of the third order term in the saddle-point expansion, corrections come from the quartic term. With $\gamma$ playing the role of the saddle-point variable, we expect an error that scales like $\frac{1}{\gamma}$, or correspondingly $\frac{1}{\delta}$.}

\rev{To explore this in more detail, we generated a range of examples by sampling random graphs of size 6, 8, 10, and 12, with density $\frac{3}{2}$ that were prime and connected. From each graph we generated a precision matrix from the $\G$-Wishart distribution \cite{mw19} and then $2n$ observations from the multivariate Gaussian to create $D$ as the data correlation matrix. Then we compute the logarithm of the ratio of normalising constants $\log\left(\frac{\const_\G(\delta, D)}{\const_\G(\delta, I_n)}\right)$ from Conjecture 1 or Equation~(\ref{eq:conj}). We compare this to a value computed for a range of $\gamma$ with an efficient Monte Carlo integration scheme \cite{GWishart} over 100 million samples for the high accuracy needed. Scripts to generate the data and plots are made available publicly\footnote{\url{https://github.com/jackkuipers/GWapprox}}.

We plot the error (Monte Carlo value minus the approximation) as a function of $\gamma$ in Figure~\ref{fig:convergeapprox} in a log-log plot. We observe very good agreement with a linear behaviour with a slope of -1, indicating that the error indeed scales like $\frac{1}{\gamma}$.}  

\rev{For each realisation of random graph and $D$, we shifted $\gamma$ to compensate for higher-order terms, which otherwise induce curvature in the log–log plots. This improves linearity and makes the leading $\gamma^{-1}$ scaling more apparent, and we optimised the linearity of each fit recording the optimal shift. The average value was around $-1$, $-2.4$, $-3.4$, and $-4.3$, for the different $n$, as displayed in Figure~\ref{fig:convergeapprox}. The average slope for the optimal fits was $-1.00$, $-1.01$, $-1.01$, $-1.03$ across the different $n$ with a standard deviations of just $0.011$, $0.013$, $0.013$, and $0.063$, respectively.}

\rev{Although every realisation generates a different line in Figure~\ref{fig:convergeapprox}, they are quite tightly packed, suggesting that the graph $\G$ and matrix $D$ have only moderate influence on the approximation error.}

\section{The next-order approximation}

\begin{figure}
	\includegraphics[width=0.48\textwidth]{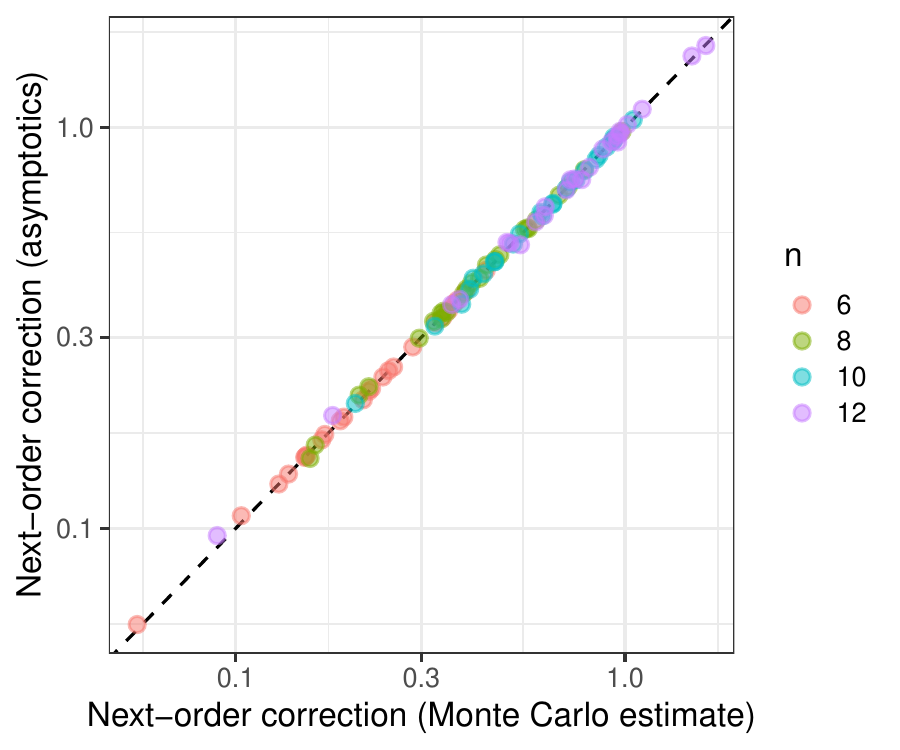}
		\caption{\rev{For the simulations depicted in Figure~\ref{fig:convergeapprox}, we extract the $\gamma^{-1}$ correction to the error ($x$-axis) and plot it against the $\left[J_{\G}(D^{\G}) - J_{\G}(I_n)\right]$ value computed from our asymptotic expansion ($y$-axis).}} \label{fig:nlo} 
\end{figure}

\begin{figure*}
\begin{tabular}{cc}
\includegraphics[width=0.48\textwidth]{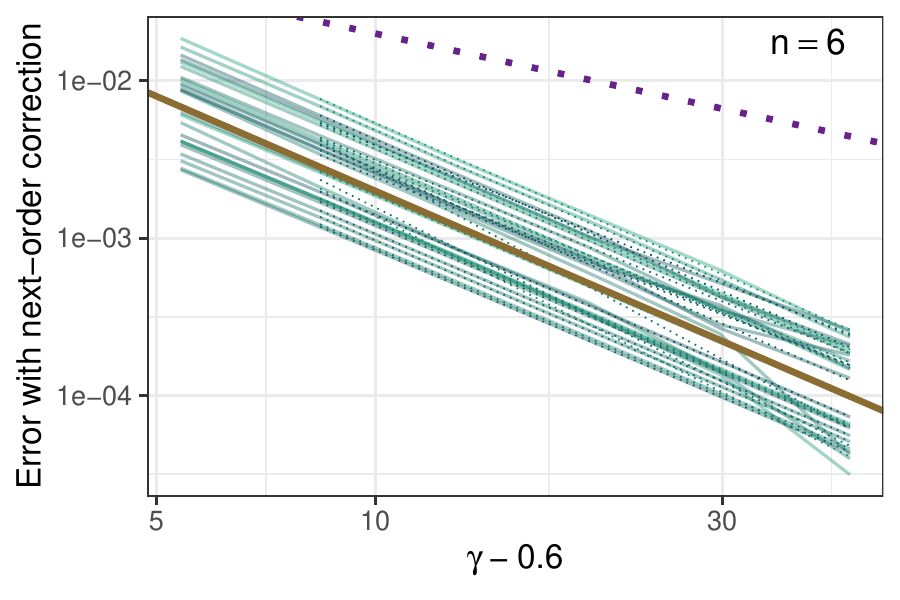} & \includegraphics[width=0.48\textwidth]{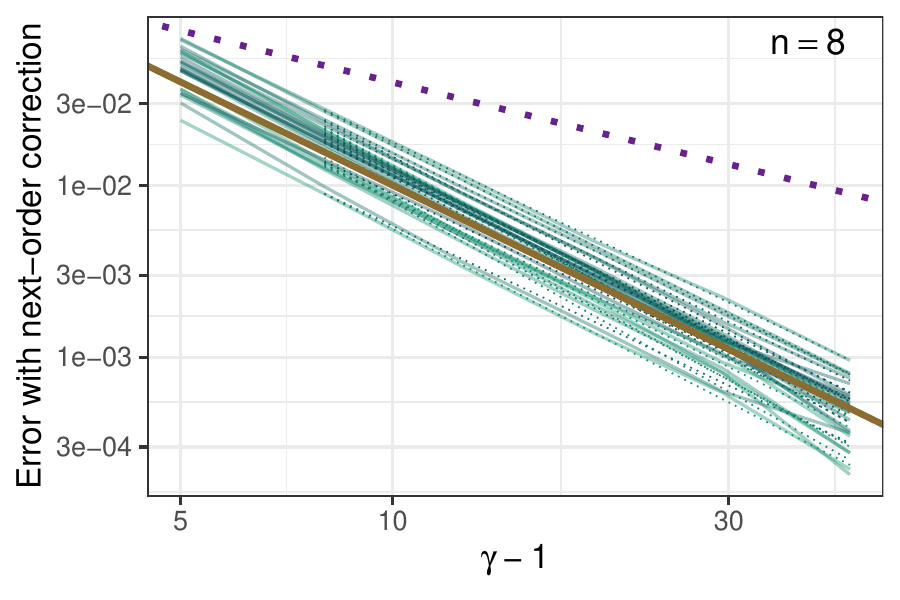}\\
\includegraphics[width=0.48\textwidth]{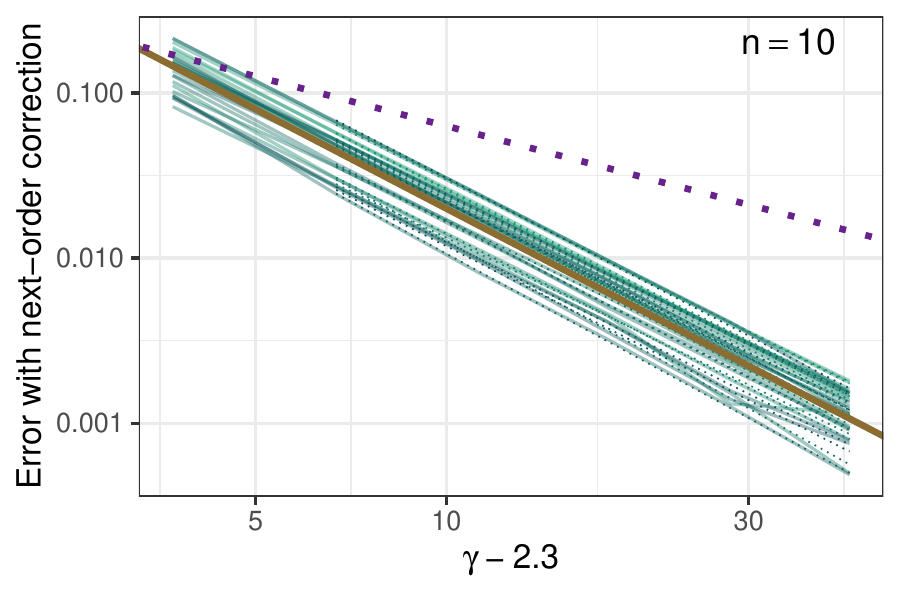} & \includegraphics[width=0.48\textwidth]{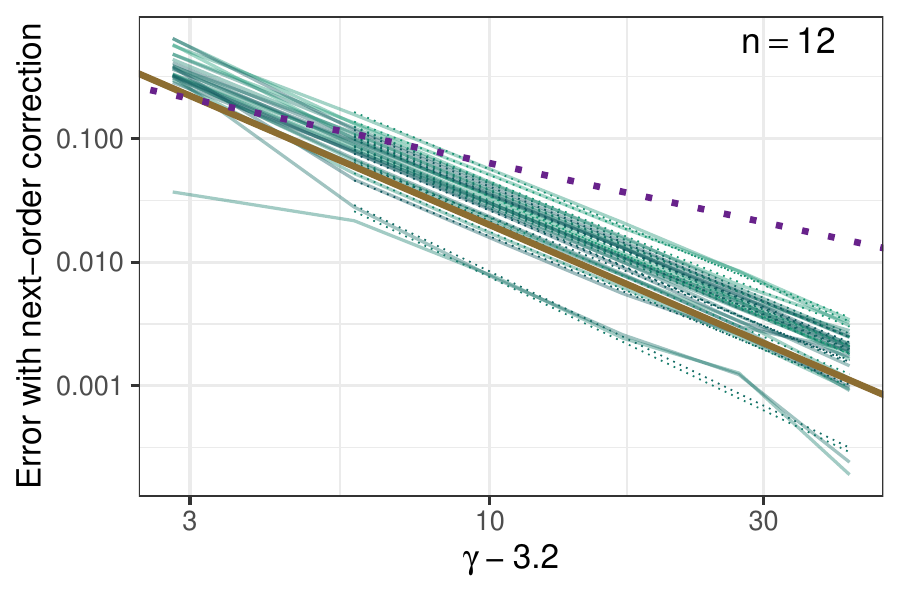}
\end{tabular}
\caption{\rev{Akin to Figure~\ref{fig:convergeapprox}, but showing the relative error after including the next-order correction of Equation~(\ref{eq:better_approx}). The indicative thick brown line has  a slope of $-2$, while we shifted $\gamma$ by $-0.6$, $-1$, $-2.3$, and $-3.2$, for the different $n$. For reference, the indicative lines of Figure~\ref{fig:convergeapprox} with slopes of $-1$ are displayed as dotted lines.}} \label{fig:convergebetterapprox} 
\end{figure*}

\rev{We turn to evaluating the next term in a $\gamma^{-1}$ asymptotic expansion and return to the central integral in Equation~(\ref{eq:int_IDG})
\begin{align}
\int_{\mathbb S^{n}(\G^c)} \det(I_n + i(D^{\G})^{-1}T)^{-\gamma} \dd T . \nonumber
\end{align}
Keeping track of higher-order terms in the expansion of the determinant we have
\begin{align}
\approx \int e^{i\frac{\gamma}{3}\tr((D^{\G})^{-1}T)^{3} + \frac{\gamma}{4}\tr((D^{\G})^{-1}T)^{4}} e^{-\gamma \boldsymbol{t}^\T \Sigma^{-1} \boldsymbol{t}} \dd \boldsymbol{t} , \nonumber
\end{align}
where the second-order term is the Gaussian as previously. Then we follow the approach of Gaussian diagrammatics \cite{novaes2024gaussian} and expand the exponential with the traces as a Taylor series 
\begin{align} \label{eq:trace_expansion}
\approx & \int\left(1 + \frac{\gamma}{4}\tr((D^{\G})^{-1}T)^{4} - \frac{\gamma^2}{2!3^2}\left[\tr((D^{\G})^{-1}T)^{3}\right]^{2}\right) \\
& \times e^{-\gamma \boldsymbol{t}^\T \Sigma^{-1} \boldsymbol{t}} \dd \boldsymbol{t}. \nonumber
\end{align}
The traces are closed paths through the matrices, alternately collecting entries of $\boldsymbol{t}$ as they pass through $T$ and weights from the $(D^{\G})^{-1}$. The $\tr((D^{\G})^{-1}T)^{4}$ is quartic in the elements of $\boldsymbol{t}$, while $\left[\tr((D^{\G})^{-1}T)^{3}\right]^{2}$ is sextic (made up of pairs of cubic terms).

We perform the Gaussian integrals over the paths from these traces using Isserlis's theorem \cite{isserlis1918formula}, also known \cite{novaes2024gaussian} as Wick's rule, which involves all possible pairing of the elements of $\boldsymbol{t}$ collected along each path and the product of the covariances of those pairs. Since the covariance is $\Sigma/2\gamma$, the trace terms in Equation~(\ref{eq:trace_expansion}) with two and three covariance terms respectively are both order $\frac{1}{\gamma}$.

Interpreting the 1 in Equation~(\ref{eq:trace_expansion}) as the leading-order contribution, the relative next-order correction is then
\begin{align} 
J_{\G}(D^{\G}) = & \frac{1}{16}\int\left(\tr((D^{\G})^{-1}T)^{4} - \frac{1}{9}\left[\tr((D^{\G})^{-1}T)^{3}\right]^{2}\right) \nonumber \\
& \times e^{-\frac{1}{2} \boldsymbol{t}^\T \Sigma^{-1} \boldsymbol{t}} \dd \boldsymbol{t} , \nonumber
\end{align}
with $\Sigma$ defined as in Equation~(\ref{eq:sigma_def}) depending on $\G$ and $D^{\G}$. Similarly we obtain the relative correction for the identity matrix,
giving us the next term in the asymptotic expansion of the ratio of normalising constants:
\begin{align} \label{eq:better_approx}
\log\left(\frac{\const_\G(\delta, D)}{\const_\G(\delta, I_n)}\right) = & 
-\gamma\log\det(D^\G) + \frac{1}{2}\log\det(\Sigma) \nonumber \\
& + \frac{1}{\gamma}\left[J_{\G}(D^{\G}) - J_{\G}(I_n)\right] + O(\gamma^{-2}).
\end{align} 

In terms of computation, for $J_{\G}(D^{\G})$ we implement a simple tracking of all paths for the quartic term, though we simplify the Wick's contractions for the sextic term since it is the square of cubic terms. For $J_{\G}(I_n)$, the quartic term just counts cycles of length 4 which cover each edge twice. These either follow the same edge back and forth, or connect through common neighbours of any pair of vertices, and are easy to count. The cubic terms correspond to triangles among the edges, and the sextic term is also straightforward to compute. We integrate code for computing the correction into the \textit{GWnorm}\footnote{\url{https://CRAN.R-project.org/package=GWnorm}} R package.

From the previous simulation results, we extract the $\gamma^{-1}$ correction by fitting $\gamma$ to the inverse error (the error is depicted in Figure~\ref{fig:convergeapprox}) and extracting the regression coefficient (the offset absorbs higher-order effects). Plotting these values against our $\left[J_{\G}(D^{\G}) - J_{\G}(I_n)\right]$ term (Figure~\ref{fig:nlo}) we observe near perfect agreement.

In fact, deviations are predominantly due to the difficulty in estimating the normalising constants to high accuracy for larger $\gamma$. Even with the large number of Monte Carlo samples employed, and using the more efficient and accurate Monte Carlo scheme \cite{GWishart}, the small stochastic noise starts to approach the magnitude of the $\gamma^{-2}$ approximation error of Equation~(\ref{eq:better_approx}) for $\gamma \approx 100$.

To visualise this, we repeat Figure~\ref{fig:convergeapprox}, but now depicting the error after using our more accurate approximation of Equation~(\ref{eq:better_approx}) in Figure~\ref{fig:convergebetterapprox}. The slopes are around $-2$, aligning with a remaining $\gamma^{-2}$ error, with average values of $-2.00$, $-2.09$, $-1.96$, $-1.96$ for the different $n$ respectively. Deviations from linear fits for larger $\gamma$ at smaller errors are also visible, indicating the limits of accuracy of the Monte Carlo scheme. For $\gamma \gtrsim 30$, our next-order approximation already gives high accuracy, and for larger $\gamma$ is increasingly favourable over Monte Carlo \cite{Atay-Kayis, GWishart} alternatives.

}

\section{Conclusions}

While efficient to compute for $\G$-Wishart normalising constants, Roverato's conjecture does not hold for all graphs. Instead we showed that the conjecture implies an approximation that had previously been employed in Bayesian samplers to speed up inference when computing the ratio of normalising constants with identity scale matrices \cite{mml23}. In Bayesian inference the parameter of the normalising constant \rev{$\delta$} is increased by the amount of data, so that although not exact, Roverato's conjecture may still offer a good approximation for general data matrices.
\rev{In fact, building on our previous Fourier analysis \cite{GWishart}, we could show that the conjecture corresponds to a saddle-point approximation, and hence captures the leading order term in an asymptotic expansion.

Connecting to Gaussian diagrammatics \cite{novaes2024gaussian} of random matrix theory, we could then derive the next-order correction. This agrees well with values estimated from Monte Carlo integration \cite{GWishart}. Obtaining accurate enough estimates from sampling as we increase the asymptotic parameter $\delta$, and increase the accuracy of the asymptotic expansion, becomes increasingly challenging. For even moderate values of $\delta$, relevant for Bayesian inference, our improved approximation may well suffice in practice.

The next-order correction of $\const_\G(\delta, D)$ depends intricately on the graph $\G$ and the data matrix $D$, and does not seem expressible in terms of $\const_\G(\delta, I_{|\V(\G)|})$. For its calculation, we need to follow paths through the complement of the graph (and the data matrix) involving at least four edges. Denoting by $\tau$ the number of edges in the complement, this suggests a complexity of $O(\tau^4)$. Compared to the complexity $O(\tau^3)$ of computing the determinant of the Isserlis matrix of Equation~(\ref{eq:conj}) as part of the leading-order terms, suggests an increase of $O(\tau)$. Higher-order corrections would involve longer paths, with corresponding increases in complexity. As $\tau$ can be $O(n^2)$, then computing next and higher-order corrections can suffer from high complexity. Creating accurate approximations, while keeping the computational cost low is then a key challenge.

For the leading-order terms, we can choose between the Isserlis matrix of the graph of Conjecture 1, and the dual version of the complement of Equation~(\ref{eq:conj}), with the former being lower complexity for sparse graphs. It may be interesting to see if the next-order correction can similarly be expressed in terms of the edges of the graph instead.

Finally, the conjecture and the approximations we developed are for the ratio of normalising constants for the same graph $\G$ with different data matrices. For Bayesian samplers of Gaussian graphical models the ratio of normalising constants with the same matrix $D$ for different graphs (akin to the approximation of \cite{mml23}) is often more pertinent. The approaches explored in this manuscript may generalise well to those tasks also.
}

\begin{acks}[Acknowledgments]
\rev{The authors thank Marcel Novaes for insightful discussions on Gaussian diagrammatics.}
\end{acks}
\begin{funding}
The first author is grateful for funding support for this work from the University of Basel through the Research Fund for Excellent Junior Researchers.
\end{funding}

\appendix

\rev{\section{Fourier analysis of the normalising constant integral}

Starting from Equation~(\ref{eq:C_G}), we first make a change of variables $K \to 2K$
\begin{equation} \nonumber
	\const_\G(\delta, D) := 2^{\frac{n\delta}{2}+\vert \E(\G) \vert}\int_{\mathbb S^n_{++}(\G)} \det(K)^{\beta} e^{-\tr(KD)} \dd K,
\end{equation}
with $\beta = \frac{\delta-2}{2}$. We are integrating over matrices which are zero when there is no edge in $\G$ reflecting the constraint for Gaussian graphical models that the precision matrix is zero for unconnected nodes. The key idea is to enforce these zeros using the matrix Dirac delta
\begin{align} \label{eq:matrixDdelta}
\delta_{\G^c}(X) = \pi^{-\vert \E(\G^c) \vert} \int_{\mathbb S^{n}(\G^c)}  e^{i \tr (TX)} \dd T ,
\end{align}
so that we can lift the restriction from the integral and integrate freely over positive definite matrices
\begin{align} \nonumber
	\const_\G(\delta, D) := 2^{\frac{n\delta}{2}+\vert \E(\G) \vert}&\int_{\mathbb S^n_{++}} \delta_{\G^c}(K) \det(K)^{\beta} \\
& \times e^{-\tr(KD)} \dd K. \nonumber
\end{align}
Substituting for the matrix Dirac delta from Equation~(\ref{eq:matrixDdelta}) and exchanging the order of integration we have
\begin{align} \nonumber \label{eq:C_Ghalfway}
\const_\G(\delta, D) := \frac{2^{\frac{n\delta}{2}+\vert \E(\G) \vert}}{\pi^{\vert \E(\G^c) \vert}}&\int_{\mathbb S^{n}(\G^c)}\int_{\mathbb S^n_{++}} \det(K)^{\beta}  \\
& \times e^{-\tr(KD)} e^{i \tr (TK)}\dd K \dd T. 
\end{align}
The central integral is just the Wishart integral over a complete graph 
\begin{align} \nonumber
&\int_{\mathbb S^n_{++}} \det(K)^{\beta} e^{-\tr(K[D-iT])}\dd K \\
& = \Gamma_n[\gamma] \det(D-iT)^{-\gamma} ,\nonumber
\end{align}
where $\gamma = \beta + \frac{n+1}{2} = \frac{n+\delta-1}{2}$ and $\Gamma_n$ is the multivariate gamma function. Substituting into Equation~(\ref{eq:C_Ghalfway})
\begin{align} \nonumber
\const_\G(\delta, D) := &\frac{2^{\frac{n\delta}{2}+\vert \E(\G) \vert}}{\pi^{\vert \E(\G^c) \vert}}\Gamma_n[\gamma] \\
&\times\int_{\mathbb S^{n}(\G^c)} \det(D-iT)^{-\gamma} \dd T,\nonumber
\end{align}
we recover Equation~(\ref{eq:I_G}) after conjugating (as the integral is real).}



\bibliographystyle{imsart-number} 
\bibliography{bibliography}       


\end{document}